\documentclass{scrartcl}

\usepackage[utf8]{inputenc} 
\usepackage[T1]{fontenc}  
\usepackage[english]{babel}
\usepackage{amsmath} 
\usepackage{amsfonts} 
\usepackage{amssymb}
\usepackage{mathrsfs}
\usepackage{authblk}
\usepackage{cite}
\usepackage{typearea}                
\usepackage{todonotes}
\usepackage{scrlayer-scrpage} 
\usepackage{tikz,pgffor}
\usetikzlibrary{shadows}
\usepackage{mathtools}
\usepackage{tabularx}
\usepackage{rotating,graphics, graphicx, wrapfig}
\usepackage{pgffor}
\usepackage{a4wide}
\usepackage{graphicx}
\graphicspath{ {figures/} }
\usepackage{makeidx}
\usepackage[ruled]{algorithm2e} 
\usepackage{ gensymb }
\usepackage[makeroom]{cancel}
\makeindex

\ohead{\headmark}
\automark[subsection]{section}

\usetikzlibrary{decorations.pathreplacing}
\usetikzlibrary{topaths,calc}
\usepackage{float}
\usetikzlibrary{snakes,arrows,patterns}

\newtheorem{defsatzusw}{}[section]

\newtheorem{theorem}[defsatzusw]{Theorem}

\newtheorem{lemma}[defsatzusw]{Lemma}
\newtheorem{corollary}[defsatzusw]{Corollary}

\newtheorem{conjecture}[defsatzusw]{Conjecture}

\newcommand{\N}{\mathbb N}

\newenvironment{proof}{\paragraph{$\mathit{Proof.}$}
	\vspace{-.4cm}\hspace{-.3cm}}{\vspace{-.2cm}\begin{flushright}$\Box$\end{flushright}}

\newcommand{\cupdot}{\mathbin{\mathaccent\cdot\cup}}

\definecolor{mygray}{HTML}{E6E6E6}
\definecolor{mygreen}{HTML}{298A08}
\definecolor{myred}{HTML}{B40404}

\pgfdeclarepatternformonly{newcrosshatch}{\pgfqpoint{-1pt}{-1pt}}{\pgfqpoint{4pt}{4pt}}{\pgfqpoint{3pt}{3pt}}%
{
	\pgfsetstrokeopacity{0.7}
	\pgfsetlinewidth{0.3pt}
	\pgfpathmoveto{\pgfqpoint{3.2pt}{0pt}}
	\pgfpathlineto{\pgfqpoint{0pt}{3.2pt}}
	\pgfpathmoveto{\pgfqpoint{0pt}{0pt}}
	\pgfpathlineto{\pgfqpoint{3.2pt}{3.2pt}}
	\pgfusepath{stroke}
}

\title{Counterexamples to the characterisation of graphs with equal independence and annihilation number}
\author[]{Michaela Hiller}
\affil[]{\small Lehrstuhl II f\"ur Mathematik, RWTH Aachen, 52062 Aachen, Germany }
\date{\vspace{-5ex}}

\begin{document}
	\setkomafont{sectioning}{\normalcolor\bfseries}	
	\maketitle

	\begin{abstract}\noindent
		\textbf{Abstract.}
		We disprove the characterisation of graphs with equal independence and annihilation number by Larson and Pepper \cite{LarsonPepper2011}. 
		Series of counterexamples with arbitrary number of vertices, arbitrary number of components, arbitrary large independence number and arbitrary large difference between the critical and the regular independence number are provided. 
		Furthermore, we point out the error in the proof of the theorem.
		However, we show that the theorem still holds for bipartite graphs and connected claw-free graphs.
	\end{abstract}

	
	\section{Introduction}
	
	In \cite{LarsonPepper2011} the authors claim the following characterisation of graphs with equal independence number $\alpha$ and annihilation number $a$ using the critical independence number $\alpha'$.
	
	\begin{theorem}\textnormal{\cite{LarsonPepper2011}}\label{falsch}
		Let $G=(V,E)$ be a graph on $n$ vertices. Then $$\alpha(G)=a(G) \quad\text{ if and only if, }\quad	
				\begin{array}{ll}
					(1)~ a(G)\geq \frac n2:& \alpha'(G)=a(G)\\[2mm]
					(2)~ a(G)= \frac{n-1}{2}:& \alpha'(G-v)=a(G) ~\text{ for some } v\in V.
				\end{array}$$
	\end{theorem}
	
	\noindent
	Since the critical independence number and the annihilation number can both be calculated in polynomial time, this result would yield a polynomial-time algorithm to verify whether the upper bound on the independence number is met for a graph.\\
	
	\noindent
	Note that the "if"-direction is still true.
	
	\begin{lemma}\label{if}
		Let $G=(V,E)$ be a graph on $n$ vertices. Then $$\alpha(G)=a(G) \quad\text{ if }\quad	
		\begin{array}{ll}
			(1)~ a(G)\geq \frac n2:& \alpha'(G)=a(G)\\[2mm]
			(2)~ a(G)= \frac{n-1}{2}:& \alpha'(G-v)=a(G) ~\text{ for some } v\in V.
		\end{array}$$
	\end{lemma} 
	\begin{proof}
		Since $\alpha'\leq\alpha\leq a$ for all graphs, (1) directly implies that $\alpha=a$. 
		In (2) we have $a(G)=\alpha'(G-v)\leq\alpha(G-v)\leq\alpha(G)\leq a(G)$. 
		Thus, all inequalities hold with equality and it follows that $\alpha=a$.\vspace{-.45cm}
	\end{proof}

%
%
	
	\noindent
	We disprove Theorem \ref{falsch} by creating various series of counterexamples in Section 2 and point out the error in the proof in Section 3. However, in Section 4 we will show that the theorem holds for restricted graph classes.

	
	\section{Counterexamples}
\begin{figure}[H]
		\begin{center}
			\begin{tikzpicture}[decoration=brace, node distance=5em, every node/.style={scale=.85}, scale=.85]
				\node(1) at (0,0)[draw, circle, scale=.85]{};
				\node(2) at (2,0)[draw, circle, scale=.85] {};
				\node(3) at (1,1.65)[circle,fill= black,scale=.85] {};
				\node(4) at (3,1.65)[circle,fill= black,scale=.85] {};
				\draw[line width=.8pt, black](1)to(2);
				\draw[line width=.8pt, black](1)to(3);
				\draw[line width=.8pt, black](3)to(2);
				
				\node(5) at (6,1.65)[]{\Large$(\bcancel{2},\bcancel{2},\cancel{2},0)$};
				\node(6) at (6.3,1.05)[]{\Large$0$};
				
				
			\end{tikzpicture}
			\caption{$C_3$ and a singleton with $\alpha=a=2\geq\frac n2$ while $\alpha'=1$.}		
			\label{C_3+Singleton}
		\end{center}
	\end{figure}
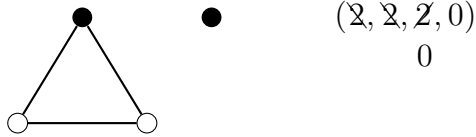

	In the following, we provide series of counterexamples with arbitrary number of vertices, arbitrary number of components, arbitrary large independence number and arbitrary large difference between the critical and the regular independence number. 
	The smallest counter example we found is a $C_3$ with a Singleton as shown in Figure \ref{C_3+Singleton}.
	We mark a maximum independence set by filled vertices.
	
	\noindent
	We can now generate counterexamples with an arbitrary number of components by adding further singletons: 
	For the graph consisting of a $C_3$ and $t$ singletons, i. e. $n=t+3$, we obtain $\alpha=a=t+1\geq\frac n2$ since the singletons yield additional vertices in each maximum independent set as well as additional zeros in the annihilation process. 
	However, $\alpha'=t$.
	
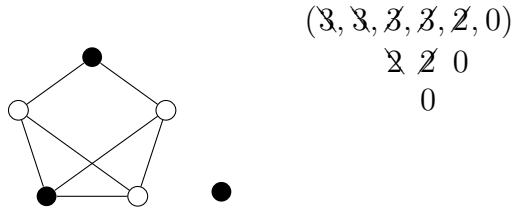
\begin{figure}[H]
		\begin{center}
			\begin{tikzpicture}[decoration=brace, node distance=5em, every node/.style={scale=.85}, scale=.85]
				\begin{scope}[shift={(3,-.1)},rotate = -3*18]
					\foreach \x [count=\p from 0] in {1,3} {
						\node(\p) [circle,fill= black,scale=.85]  at (-\x*360/5:1.2) {};};
					\foreach \x [count=\p from 2] in {0,2,4} {
						\node(\p) [draw,circle,scale=.85] at (-\x*360/5:1.2) {};};
					\draw[] (2) -- (0) -- (3) -- (1) -- (4) -- (2);
					\draw[] (2) -- (3);
					\draw[] (0) -- (4);
					
				\end{scope}
				\node(6) at (5,-1)[circle,fill= black, scale=.85] {};
				\node(5) at (7.9,1.65)[]{\Large$(\bcancel{3},\bcancel{3},\cancel{3},\cancel{3},\cancel{2},0)$};
				\node(6) at (8.18,1.05)[]{\Large$\bcancel{2}$\,\;$\cancel{2}$\:\;$0$};
				\node(6) at (8.2,.45)[]{\Large\phantom{$\bcancel{2}$}$0$\phantom{$\bcancel{2}$}};
			\end{tikzpicture}
			\caption{For a $C_5$ with two chords and a singleton, we obtain $\alpha=a=3\geq\frac n2$ but $\alpha'=1$.}
			\label{C_5+Chords+Singleton}	
		\end{center}
	\end{figure}

	A further counter example is the graph consisting of a $C_5$ with two chords and a singleton as in Figure \ref{C_5+Chords+Singleton}.
	
	\noindent
	Counterexamples are not required to contain singletons: 
	The degree sequence of an odd cycle $C_{2k+1}$, $k\in\N$ combined with an path of odd length $P_{2l+1}$, $l\in\N$ is $\left(2^{2(k+l-1)},1^2\right)$. Thus, we have $\alpha=a=4\geq\frac n2$ as $n=k+l$, but $\alpha'=l+1$.
	
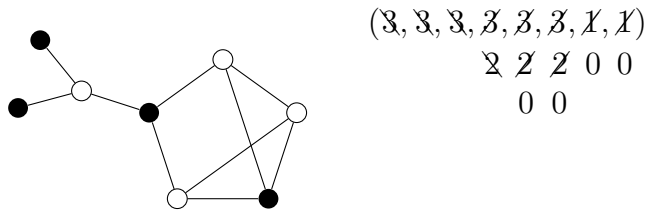
\begin{figure}[H]
		\begin{center}
			\begin{tikzpicture}[decoration=brace, node distance=5em, every node/.style={scale=.85}, scale=.85]
				\begin{scope}[shift={(3,-.1)},rotate = 18]
				\foreach \x [count=\p from 0] in {1,3} {
					\node(\p) [circle,fill= black,scale=.85]  at (-\x*360/5:1.2) {};};
				\foreach \x [count=\p from 2] in {0,2,4} {
					\node(\p) [draw,circle,scale=.85] at (-\x*360/5:1.2) {};};
				\draw[] (2) -- (0) -- (3) -- (1) -- (4) -- (2);
				\draw[] (2) -- (3);
				\draw[] (0) -- (4);
				\node(5) at (-3*360/5:2.3)[draw, circle, scale=.85] {};
				\node(6) at (-3*360/5+10:3.2)[circle,fill= black, scale=.85] {};
				\node(7) at (-3*360/5-10:3.2)[circle,fill= black, scale=.85] {};
				\draw (1) -- (5);
				\draw (6) -- (5);
				\draw (7) -- (5);
				\end{scope}
			
				\node(5) at (7.4,1.65)[]{\Large$(\bcancel{3},\bcancel{3},\bcancel{3},\cancel{3},\cancel{3},\cancel{3},\cancel{1},\cancel{1})$};
				\node(6) at (8.18,1.05)[]{\Large$\bcancel{2}$\,\;$\cancel{2}$\:\;$\cancel{2}$\:\;$0$\:\;$0$};
				\node(6) at (8.18,.45)[]{\Large\phantom{$\bcancel{2}$}\,\;$0$\;\;$0$\:\;\phantom{$0$}\:\;\phantom{$0$}};
			\end{tikzpicture}
		\caption{This graph with $n=8$ fulfils $\alpha=a=4\geq\frac n2$, while $\alpha'=2$.}
		\label{C_5+Chrods+Star}		
	\end{center}
	\end{figure}

	Furthermore, we can provide counterexamples that are connected as the graph shown in Figure \ref{C_5+Chrods+Star}.
	
	\noindent
	To obtain counterexamples where the difference between the critical independence number and the annihilation number becomes arbitrary large, the example above can be generalised: 
	Starting with a $C_{2k+1}=\{v_1,\dots,v_{2k+1}\}$ for $k\geq2$, we add chords $\{v_i,v_{i+k}\}$ for $i\in\{1,\dots,k\}$ and attach a $P_3$ with its central vertex to $v_{2k+1}$. 
	By construction, all vertices have degree three except the two end vertices of the $P_3$, which we denote by $x_1$ and $x_2$. 
	Thus, the degree sequence of this connected graph with $n=2k+4$ is $(3^{2k+2},1,1)$ yielding $a=k+2$.
	
	\begin{figure}[H]
		\begin{center}
			\begin{tikzpicture}[decoration=brace, node distance=5em, every node/.style={scale=.85}, scale=.85]
				
				\node(5) at (-.02,1.65+.65)[]{\Large$(\overbrace{\bcancel{3},\dots,\bcancel{3}}^{k},\bcancel{3},\cancel{3},\overbrace{\cancel{3},\cancel{3},\dots,\cancel{3}}^{k},\cancel{1},\cancel{1})$};
				\node(5) at (0,1.05+.3)[]{\Large$\phantom{(\bcancel{3},\dots,\bcancel{3},\bcancel{3},}\bcancel{2}\;\,\cancel{2}~\ldots\;\,\cancel{2}\;\,\cancel{2}\;\:0\;\:0\phantom{)}$};
				\node(5) at (-.07,0.45-.05)[]{\Large$\phantom{(\bcancel{3},\dots,\bcancel{3},\bcancel{3},\cancel{3}}\underbrace{0\,\,\ldots\;\,0\;\:0\phantom{\cancel{1},\cancel{1})}}_{k+2}$};
			\end{tikzpicture}
		\end{center}
	\end{figure}
\vspace{-.3cm}
\noindent	
	Clearly, any maximum independent set contains $x_1$, $x_2$ and at most $k$ vertices on the cycle, i. e. $\alpha\leq k+2$.
	It remains to prove $\alpha\geq k+2$.
	We consider cases $k\equiv_20$ and $k\equiv_21$.
	
	\noindent
	For the former case, we define an independent set
		$$I=\{x_1,x_2,v_{2k+1}, \underbrace{v_2,v_4,\dots,v_k}_{\frac k2}, \underbrace{v_{k+3},v_{k+5},\dots,v_{2k-1}}_{\frac{k-2}{2}}\}$$
	and in the latter 
		$$I=\{x_1,x_2,v_{2k+1}, \underbrace{v_2,v_4,\dots,v_{k-1}}_{\frac{k-1}{2}}, \underbrace{v_{k+1},v_{k+3},\dots,v_{2k-2}}_{\frac{k-1}{2}}\}.$$
	
	\noindent
	In both cases $|I|=k+2\leq\alpha$, thus $\alpha=k+2=a$, whereas $\alpha'=2$. Now, $\alpha-\alpha'=k$ can become arbitrary large.
	For $k=2$, the constructed graph corresponds to the one in Figure \ref{C_5+Chrods+Star}; for $k=3$ and $k=4$, the graphs are shown in Figure \ref{C_7+Chrods+Star} and Figure \ref{C_9+Chrods+Star}, respectively.
	
	\begin{figure}[H]
		\begin{center}
			\begin{tikzpicture}[decoration=brace, node distance=5em, every node/.style={scale=.85}, scale=.85]
				\begin{scope}[shift={(3,-.1)},rotate = 141]
					\foreach \x [count=\p from 0] in {0,2,4} {
						\node(\x) [circle,fill= black,scale=.85]  at (\x*360/7:1.6) {};};
					\foreach \x [count=\p from 2] in {1,3,5,6} {
						\node(\x) [draw,circle,scale=.85] at (\x*360/7:1.6) {};};
					\draw[] (0) -- (1) -- (2) -- (3) -- (4) -- (5) -- (6) -- (0);
					\draw[] (1) -- (4);
					\draw[] (2) -- (5);
					\draw[] (3) -- (6);
					\node(7) at (-0*360/7:2.7)[draw, circle, scale=.85] {};
					\node(8) at (-0*360/7+10:3.7)[circle,fill= black, scale=.85] {};
					\node(9) at (-0*360/7-10:3.7)[circle,fill= black, scale=.85] {};
					\draw (0) -- (7);
					\draw (7) -- (8);
					\draw (7) -- (9);
				\end{scope}
				
				\node(5) at (8.2,1.65+.3)[]{\Large$(\bcancel{3},\bcancel{3},\bcancel{3},\bcancel{3},\cancel{3},\cancel{3},\cancel{3},\cancel{3},\cancel{1},\cancel{1})$};
				\node(5) at (8.2,1.05+.3)[]{\Large$\phantom{(\bcancel{3},\bcancel{3},\bcancel{3},\bcancel{3},}\bcancel{2}\;\,\cancel{2}\;\,\cancel{2}\;\,\cancel{2}\;\:0\;\:0\phantom{)}$};
				\node(5) at (8.2,0.45+.3)[]{\Large$\phantom{(\bcancel{3},\bcancel{3},\bcancel{3},\bcancel{3},\cancel{3},}0\;\;0\;\;0\phantom{,\cancel{1},\cancel{1})}$};
			\end{tikzpicture}
			\caption{For $k=3$, we get $n=10$ and $\alpha=a=5\geq\frac n2$, while $\alpha'=2$.}
			\label{C_7+Chrods+Star}			
		\end{center}
	\end{figure}	
			
	\begin{figure}[H]
		\begin{center}
			\begin{tikzpicture}[decoration=brace, node distance=5em, every node/.style={scale=.85}, scale=.85]
				\begin{scope}[shift={(3,-.1)},rotate = 130]
					\foreach \x [count=\p from 0] in {0,2,4,7} {
						\node(\x) [circle,fill= black,scale=.85]  at (\x*360/9:1.9) {};};
					\foreach \x [count=\p from 2] in {1,3,5,6,8} {
						\node(\x) [draw,circle,scale=.85] at (\x*360/9:1.9) {};};
					\draw[] (0) -- (1) -- (2) -- (3) -- (4) -- (5) -- (6) -- (7) -- (8) -- (0);
					\draw[] (1) -- (5);
					\draw[] (2) -- (6);
					\draw[] (3) -- (7);
					\draw[] (4) -- (8);
					\node(9) at (-0*360/9:3.2)[draw, circle, scale=.85] {};
					\node(10) at (-0*360/9+10:4.2)[circle,fill= black, scale=.85] {};
					\node(11) at (-0*360/9-10:4.2)[circle,fill= black, scale=.85] {};
					\draw (0) -- (9);
					\draw (9) -- (10);
					\draw (9) -- (11);
				\end{scope}
				
				\node(5) at (8.9,1.95+.6)[]{\Large$(\bcancel{3},\bcancel{3},\bcancel{3},\bcancel{3},\bcancel{3},\cancel{3},\cancel{3},\cancel{3},\cancel{3},\cancel{3},\cancel{1},\cancel{1})$};
				\node(5) at (8.9,1.35+.6)[]{\Large$\phantom{(\bcancel{3},\bcancel{3},\bcancel{3}\bcancel{3},,\bcancel{3},}\bcancel{2}\;\,\cancel{2}\;\,\cancel{2}\;\,\cancel{2}\;\,\cancel{2}\;\:0\;\:0\phantom{)}$};
				\node(5) at (8.91,0.75+.6)[]{\Large$\phantom{(\bcancel{3},\bcancel{3},\bcancel{3},\bcancel{3},\bcancel{3},\cancel{3},}0\;\;0\;\;0\;\;0\phantom{,\cancel{1},\cancel{1})}$};
			\end{tikzpicture}
			\caption{For $k=4$, we get $n=12$ and $\alpha=a=6\geq\frac n2$, while $\alpha'=2$.}		
			\label{C_9+Chrods+Star}	
		\end{center}
	\end{figure}

	
	\section{Error in the proof}

	The error in the proof in \cite{LarsonPepper2011} occurs in the case, where $G$ is not empty, $a(G)\geq\frac{n}{2}$, the neighbourhood of the maximum critical independent set $J$ is empty and $a(G-J)<\frac{n(G-J)}{2}$. 
	The authors use the inductive assumption on $G-J+u$ for a vertex $u\in J$, but for $J=\{u\}$, i. e. $|J|=1$, we have $G-J+u=G$ and the inductive assumption cannot be applied. 
	Since the theorem is proven by induction, it is unclear whether the given proof is salvageable even for restricted graph classes.
	Since the Theorem \ref{falsch} is not true in general, the proof of the corollary (Theorem 3.3 in \cite{LarsonPepper2011}), in which König-Egerváry graphs, i. e. graphs with $\alpha+\mu=n$, are characterised by the equality of independence and annihilation number, is invalid as well. 
	
	\begin{corollary}\textnormal{\cite{LarsonPepper2011}}
		For a graph $G$ with $a(G)\geq\frac n2$, $\alpha(G)=a(G)$ if an only if $G$ is a König-Egerváry graph and every maximum independent set of $G$ is a maximum annihilating set.
	\end{corollary}

	\noindent
	Consider for example the graphs constructed in Section 2 (see Figure \ref{C_5+Chrods+Star}, \ref{C_7+Chrods+Star}, \ref{C_9+Chrods+Star}) with $n=2k+4$ and $\alpha=a=k+2\geq \frac n2$. 
	The matching number of such a graph is $\mu=k+1$.
	Thus $\alpha+\mu<n$ and the graph is not König-Egerváry.
	This also implies that the "only if"-part of Conjecture 3.4 in \cite{LarsonPepper2011} is not true.
	
	\begin{conjecture}\textnormal{\cite{LarsonPepper2011}}
		For a graph $G$ with $a\geq\frac n2$, $\alpha=a$ if and only if $G$ is a König-Egerváry graph and every maximum independent set of $G$ is a maximal annihilating set.
	\end{conjecture}

	\noindent
	In \cite{LevitMandrescu2018} and \cite{LevitMandrescu2020}, the authors gave counterexamples for the "if"-direction, but showed the "only if"-direction using the above disproved  results by Larson and Pepper. 

	It remains to retrace in which papers the theorems have been used beyond the above mentioned and to review whether the subsequent results still hold.

\newpage


	\section{Theorem for bipartite graphs and connected claw-free graphs}

	It is striking that all counterexamples mentioned above are either non-connected or contain a claw and an odd cycle.
	In fact, it turns out that the theorem still holds for bipartite graphs and connected claw-free graphs.\\
	
	\noindent
	Note that for bipartite graphs that case $a=\frac{n-1}{2}$ cannot occur. Thus, for this graph class the following theorem is equivalent to Theorem \ref{falsch}.
	
	\begin{theorem}
		Let $G$ be a bipartite graph. Then $$\alpha(G)=a(G) \quad\text{ if and only if, }\quad	
		\alpha'(G)=a(G).$$
	\end{theorem}
	
	\begin{proof}
		As seen in Lemma \ref{if} the "if"-direction holds true. 
		Hence, it remains to prove that if $G$ is bipartite and $a\geq \frac n2$, the equality $\alpha=a$ implies $\alpha'=a$.
		For bipartite graphs $\alpha'=\alpha$ \cite{Larson2011}. Therefore, the implication is true.\vspace{-.54cm}
	\end{proof}
	
	\noindent
	To prove the theorem for claw-free graphs, we need the following lemma.
	
	\begin{lemma}\label{exist}
		Let $G=(V,E)$ be connected and claw-free with $a(G)=\frac{n-1}{2}$. Then there exists a vertex $v\in V$ that does not occur in every maximum independent sets, while $G-v$ is still connected.
	\end{lemma}
	\begin{proof}
		Let $I\subseteq V$ be a maximum independent set in $G$.
		Then $V=I \cupdot \mathcal{N}(I)$.
		Consider a path $P=\{v_0,\dots,v_k\}$ of maximum length in $G$.
		If $v_0\not\in I$, the removal of $v_0$ preserves the connectedness of $G$ since all neighbours of $v_0$ have to be in $P$; otherwise, $P$ was not a path of maximum length. 
		If $v_0\in I$, then $v_1\not\in I$.
		Note that $a=\frac{n-1}{2}$ implies $\sum_{v\in X}\deg(v)<\sum_{v\in Y}\deg(v)$ for all $X,Y\subseteq V$, $X\cap Y=\emptyset$ with $|X|<|Y|$.
		Therefore, the minimum degree $\delta\geq 2$ and thus $v_0$ has at least one additional neighbour (apart from $v_1$), which as seen above has to be in $P$.
		Now assume that $G-v_1$ was not connected. 
		Then there exists a neighbour $w$ of $v_1$ that is not adjacent to any vertex of $P$.
		From $\delta\geq 2$ it follows that $w$ has another neighbour $z\not\in P$.
		But this contradicts the assumption that $P$ is a path of maximum length since we obtain a longer path by replacing $v_0$ with $w$ and $z$.\vspace{-.53cm}
	\end{proof}
	
	\begin{theorem}
		Let $G$ be a connected claw-free graph. Then $$\alpha(G)=a(G) \quad\text{ if and only if, }\quad	
		\begin{array}{ll}
			(1)~ a(G)\geq \frac n2:& \alpha'(G)=a(G)\\[2mm]
			(2)~ a(G)= \frac{n-1}{2}:& \alpha'(G-v)=a(G) ~\text{ for some } v\in V.
		\end{array}$$
	\end{theorem}
	\begin{proof}
		By Lemma \ref{if}, it suffices to consider the "only if"-direction for claw-free graphs $G=(V,E)$.\\
		First, let $\alpha(G)=a(G)\geq\frac n2$ and suppose $\alpha'(G)<a(G)$.
		Then, there exists a maximum critical independent set $I^c\subseteq V$ with $|I^c|-|\mathcal{N}(I^c)|>0$ as well as $|I^c|<|I|$ and $|I^c|-|\mathcal{N}(I^c)|>|I|-|\mathcal{N}(I)|$ for all maximum independent sets $I\subseteq V$.
		This implies $|I^c|\geq2$, $|\mathcal{N}(I^c)|\geq1$ and for $R:=V\setminus(I^c\cup\mathcal{N}(I^c))$ we get $|R|\geq3$.
		By assumption $G$ is claw-free and $I^c$ is an independent set. 
		Thus each vertex in $\mathcal{N}(I^c)$ can have at most two neighbours in $I^c$.\\
		We consider the bipartite graph $\tilde{G}:=(I^c\cup\mathcal{N}(I^c), \tilde{E})$ with $\tilde{E}:=\{uv\in E~|~u\in \mathcal{N}(I^c),v\in I^c\}$.
		Any connected component $K$ of $\tilde{G}$ falls naturally into one of two types; it is 
		$|K\cap I^c|\leq|K\cap\mathcal{N}(I^c)|$ or 
		$|K\cap I^c|>|K\cap\mathcal{N}(I^c)|$. 
		We call the former Type I components and the latter Type II components.\\
		Note that there exists at least one component of each type: 
		Since $|I^c|-|\mathcal{N}(I^c)|>0$, there is a vertex in $\mathcal{N}(I^c)$ which is adjacent to exactly two vertices in $I^c$. And as $G$ is claw-free, this vertex cannot have a neighbour in $R$.
		Additionally, since $G$ is connected, there is a vertex in $\mathcal{N}(I^c)$ which is adjacent to exactly one vertex in $I^c$ and at least one vertex in $R$.\\
		In a Type II component $K$, we have at least $|K\cap I^c|+|K\cap\mathcal{N}(I^c)|-1\geq2|K\cap\mathcal{N}(I^c)|$ edges. 
		Therefore, all vertices of $K$ in $\mathcal{N}(I^c)$ have degree two in $\tilde{G}$; otherwise, there also existed a vertex in $\mathcal{N}(I^c)$ of degree three, contradicting the claw-freeness.
		Furthermore, Type II components cannot contain vertices adjacent to $R$ in $G$; otherwise, $G$ is not claw-free. 
		Therefore, at least one Type II component $K_2$ has to be connected to a Type I component $K_1$ by an edge between two vertices in $\mathcal{N}(I^c)$.
		Let $w_1\in K_1\cap \mathcal{N}(I^c)$ and $w_2\in K_2\cap \mathcal{N}(I^c)$ be these vertices with $w_1w_2\in E$.
		But since $w_2$ has already two non-connected neighbours within $K_2$, this contradicts the claw-freeness of $G$. Hence, the claim is proven to be true for $a(G)\geq\frac n2$.\\
		Now, let $\alpha(G)=a(G)=\frac{n-1}{2}$. 
		By Lemma \ref{exist}, there exists a vertex $v\in V$ that does not occur in every maximum independent sets, while $G-v$ is still connected. 
		Thus, $\alpha(G-v)=\alpha(G)=a(G)=\frac{n-1}{2}$.
		As $\alpha(G-v)\leq a(G-v)\leq a(G)$, it follows that $\alpha(G-v)=a(G-v)=\frac{n-1}{2}$.
		Since $G$ is assumed to be claw-free, $G-v$ is claw-free as well. 
		Further, $a(G-v)\geq \frac{n(G-v)}{2}$ and therefore, we can apply the first case to $G-v$ and obtain that $\alpha'(G-v)=a(G-v)=\frac{n-1}{2}=a(G)$.
		This completes the proof for claw-free graphs.\vspace{-.43cm}
	\end{proof}
	
	\noindent
	Note that it was already shown in 1980 by Sbihi \cite{Sbihi1980} that maximum independent sets can be found in claw-free graphs in polynomial time. 
	The proof uses the blossom algorithm by Edmonds from 1965 \cite{Edmonds1965}, which yields maximum matchings in polynomial time for any graph: 
	Any maximum matching in a graph translates to a maximum independent set in the corresponding line graph and all claw-free graphs can be considered as the line graph of some graph.
	
	It remains open whether the theorem holds for other restricted graph classes or for arbitrary graphs with $a=\frac{n-1}{2}$. Of particular interest in this regard are graph classes for which maximum independent sets cannot be found in polynomial time.

	\bibliographystyle{alpha} 
	\bibliography{Literatur_Erratum}
\end{document}